\documentclass[12pt]{article}
\usepackage[letterpaper,text={6.75in,9.25in},centering]{geometry}      
\usepackage{graphicx}
\usepackage{amssymb}
\usepackage{epstopdf}
\usepackage[usenames]{color}

\usepackage{amsfonts,amsmath,amsthm}
\usepackage{caption,graphicx}
\usepackage{enumerate}
\usepackage{comment}

\newcommand{\real}{\mathbb{R}}

\newcommand{\reff}[1]{(\ref{#1})}

\newcommand{\cLm}{{ {\cal L}_{\mu}}}

\newtheorem{Theorem}{Theorem}              
\newtheorem{Lemma}{Lemma}[section]

\newtheorem{Remark}{Remark}[section]

\specialcomment{genescomments}
\specialcomment{margaretscomments}

\includecomment{genescomments}
\includecomment{margaretscomments}

\begin{document}

\title{Using global invariant manifolds to understand metastability in Burgers equation with small viscosity}
\author{Margaret Beck\thanks{email: Margaret\_Beck@brown.edu; The majority of this work was done while MB was affiliated with the Department of Mathematics, University of Surrey, Guildford, GU2 7XH, UK} \\
Division of Applied Mathematics \\
Brown University \\ 
Providence, RI~02912, USA \\[1cm]
C. Eugene Wayne\thanks{email: cew@math.bu.edu}\\
Department of Mathematics and Center for BioDynamics\\
Boston University\\
Boston, MA~02215, USA 
}                                                         

%\date{}                                           % Activate to display a given date or no date

\maketitle

\begin{abstract}
The large-time behavior of solutions to Burgers equation with small viscosity is described using invariant manifolds. In particular, a geometric explanation is provided for a phenomenon known as metastability,
which in the present context means that solutions spend a very long time near the family of solutions known as diffusive N-waves before finally converging to a stable self-similar diffusion wave. More precisely, it is shown that in terms of similarity, or scaling, variables in an algebraically weighted $L^2$ space, the self-similar diffusion waves correspond to a one-dimensional global center manifold of stationary solutions. Through each of these fixed points there exists a one-dimensional, global, attractive, invariant manifold corresponding to the diffusive N-waves. Thus, metastability corresponds to a fast transient in which solutions approach this ``metastable" manifold of diffusive N-waves, followed by a slow decay along this manifold, and, finally, convergence to the self-similar diffusion wave. 
\end{abstract}

\newpage

\section{Introduction}\label{sec_intro}

It is well known that viscosity plays an important role in the evolution of solutions to viscous conservation laws and that its presence significantly impacts the asymptotic behavior of solutions. Much work has been done to understand the relationship between solutions for zero and nonzero viscosity. 
For an overview, see, for example, \cite{Dafermos05, Liu00}. With regard to Burgers equation, one key property is the following. If $u^\mu= u^\mu(x,t)$ denotes the solution to Burgers equation with viscosity $\mu$ and $u^0= u^0(x,t)$ denotes the solution to the inviscid equation, then it is known that $u^\mu \to u^0$ in an appropriate sense for any fixed $t>0$ as $\mu \to 0$. However, for fixed $\mu$, the large time behavior of $u^\mu$ and $u^0$ is quite different, and they converge to solutions known as diffusion waves and N-waves, respectively. Thus, the limits $\mu \to 0$ and $t \to \infty$ are not interchangeable.  

Recently, a phenomenon known as metastability has been observed in Burgers equation with small viscosity on an unbounded domain \cite{KimTzavaras01}. Generally speaking, metastable behavior is when solutions exhibit long transients in which they remain close to some non-stationary state (or family of non-stationary states) for a very long time before converging to their asymptotic limit.
In \cite{KimTzavaras01}, the authors observe numerically that solutions spend a very long time near a family of solutions known as ``diffusive N-waves," before finally converging to the stable family of diffusion waves. This terminology\footnote{These diffusive N-waves are also discussed in \cite[\S4.5]{Whitham99}, where they are referred to simply as N-waves. Here, as in \cite{KimTzavaras01}, we reserve the term N-wave for solutions of the inviscid equation and diffusive N-wave for solutions of the viscous equation.} is due to the fact that the diffusive N-waves are close to inviscid N-waves. In \cite{KimTzavaras01} this is proven in a pointwise sense. Furthermore, in terms of scaling, or similarity, variables, they compute an asymptotic expansion for solutions to Burgers equation with small viscosity. They find that the stable diffusion waves correspond to the first term in the expansion, whereas the diffusive N-waves correspond to taking the first two terms. Thus, by characterizing the metastability in terms of these diffusive N-waves, they provide a way of understanding the interplay between the limits $\mu \to 0$ and $t \to \infty$.
 
\begin{figure}
\centering\includegraphics[scale=0.6,angle=0]{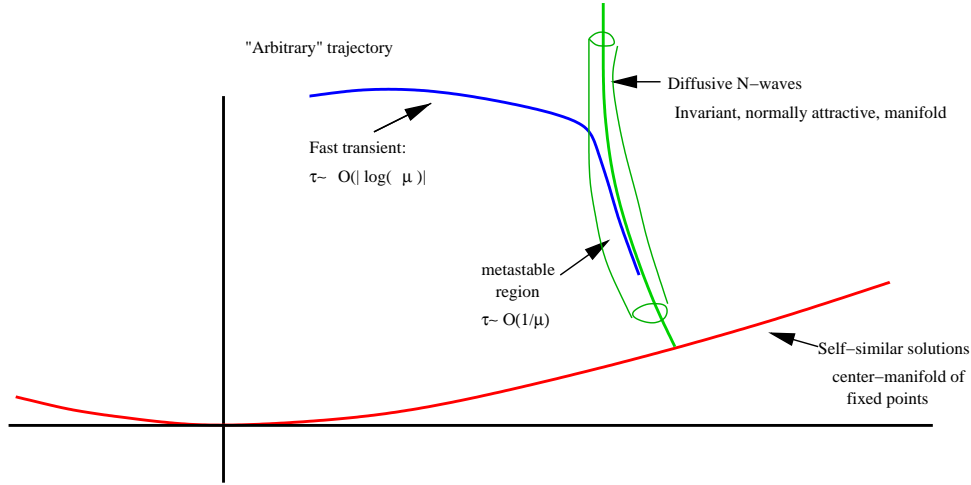}
\caption{A schematic diagram of the invariant manifolds in the phase space of Burgers equation, (\ref{eqn:scaling}), and their role in the metastable behavior. The solution trajectory (blue) experiences an initial fast transient of $\tau = \mathcal{O}(|\log \mu|)$ before entering a neighborhood of the manifold of diffusive N-waves (green). It then remains in this neighborhood for all time as it approaches, on the slower time scale of $\tau=\mathcal{O}(1/\mu)$, a point on the manifold of stable stationary states (red).}  
\label{F:geometry}
\end{figure}

In this paper, we show that the metastable behavior in the viscous Burgers equation, described in \cite{KimTzavaras01}, can be viewed as the approach to, and motion along, a normally attractive, invariant manifold in the phase space of the equation. In terms of the similarity variables, we show that one has the following picture. There exists a global, one-dimensional center manifold of stationary solutions corresponding to the self-similar diffusion waves. Through each of these fixed points there exists a global, one-dimensional, invariant, normally attractive manifold corresponding to the diffusive N-waves. For almost any initial condition, the corresponding solution of Burgers equation
approaches one of the diffusive N-wave manifolds on a relatively fast time scale: 
$\tau = \mathcal{O}(|\log\mu|)$. Due to attractivity, the solution remains close to this manifold for all time and moves along it on a slower time scale, $\tau = \mathcal{O}(1/\mu)$, towards the fixed point which has the same total mass. Note this this corresponds to an extremely long timescale $t \approx {\cal O}(e^{1/\mu})$ in
the original unscaled time variable. This scenario is illustrated in Figure \ref{F:geometry}. 

Note that, in \cite{Liu00}, it was shown that the large time behavior of solutions to a general class of conservation laws is governed by that of solutions to Burgers equation. Roughly speaking, this is due to the marginality of the nonlinearity in the case of Burgers equation and the fact that any higher order nonlinear terms in other conservation laws are irrelevant. Therefore, the present analysis for Burgers equation could potentially be used to predict and understand metastability in other conservation laws with small viscosity, as well.

We remark that a similar metastable phenomenon has also been investigated in Burgers equation on a bounded interval \cite{BerestyckiKaminSivashinsky95, SunWard99} and numerically observed in the Navier-Stokes equations on a two-dimensional bounded domain \cite{YinMontgomeryClercx03}. Furthermore, metastability has been observed in reaction-diffusion equations: for example, on a bounded interval \cite{FuscoHale89, CarrPego89, CarrPego90} and spatially discrete lattice \cite{GrantVanVleck95}. In \cite{CarrPego90}, the metastable states were described in terms of global unstable invariant manifolds of equilibria. 

The remainder of the paper is organized as follows. In \S\ref{S:set_up_results} we state the equations and function spaces within which we will work, as well as some preliminary facts about the existence of invariant manifolds in the phase space of Burgers equation. We also precisely formulate our results, Theorems \ref{thm:initial_transient} and \ref{thm:locally_attractive}. In sections \S\ref{sec_lemma}, \S\ref{sec_transient}, and \S\ref{local_attractivity} we prove Lemma \ref{lem:N_wave_closeness}, Theorem \ref{thm:initial_transient}, and Theorem \ref{thm:locally_attractive}, respectively. Concluding remarks are contained in \S\ref{S:discussion}. Finally, the appendix contains a calculation that is referred to in \S\ref{S:set_up_results}. 

%%%%%%%%%%%%%%%%%%%%%%%%%%%%%%%%%%%%%%%%%%%%%%%%%

\section{Set-up and statement of results} \label{S:set_up_results}

We now explain the set-up for the analysis and some preliminary results on invariant manifolds. Our main results are precisely stated in \S\ref{S:main_results}.

\subsection{Equations and scaling variables}

The scalar, viscous Burgers equation is the initial value problem
\begin{eqnarray}
\partial_t u &=& \mu \partial_{x}^2 u - u u_x \nonumber \\
u|_{t=0} &=& h \mbox{,} \label{burgers1}
\end{eqnarray}
and we assume the viscosity coefficient $\mu$ is small: $0< \mu \ll 1$. 
For reasons described below, it is convenient to work in so-called similarity or scaling variables, defined as
\begin{gather} 
u(x,t) = \frac{1}{\sqrt{1+t}} w\left(\frac{x}{\sqrt{1+t}},\log(t+1)\right) \nonumber  \\
\xi= \frac{x}{\sqrt{1+t}} \mbox{, \quad} \tau = \log(t+1) \mbox{.}
\label{sim}
\end{gather}
In terms of these variables, equation \reff{burgers1} becomes
\begin{eqnarray}
\partial_{\tau} w &=& \cLm w - w w_{\xi} \nonumber \\
w|_{\tau=0} &=& h \mbox{,} \label{eqn:scaling}
\end{eqnarray}
where $\cLm w = \mu \partial_\xi^2w + \frac{1}{2} \partial_{\xi} (\xi w)$.

We will study the evolution of \reff{eqn:scaling} in the algebraically weighted Hilbert space
\begin{equation*}
L^2(m) =  \left\{ w \in L^2({\mathbb{R}}) : \|w\|_{L^2(m)}^2 = \int (1+ \xi^2)^m |w(\xi)|^2 d\xi < \infty \right\} .
\label{E:def_Hs}
\end{equation*}
It was shown in \cite{GallayWayne02} that, in the spaces $L^2(m)$ with $m > 1/2$, the operator $\cLm$ generates a strongly continuous semigroup and its spectrum is given by
\begin{equation}
\sigma(\cLm) = \left\{ -\frac{n}{2}, n \in \mathbb{N} \right\} \cup \left\{\lambda \in \mathbb{C} : \mbox{Re}(\lambda) \le \frac{1}{4} - \frac{m}{2}\right\} \mbox{.} \label{E:spectrum_Lm}
\end{equation}
This is exactly the reason why the similarity variables are so useful. Equation (\ref{E:spectrum_Lm}) shows that the operator $\cLm$ has a gap, at least for $m > 1/2$, between the continuous part of the spectrum and the zero eigenvalue. As $m$ is increased, more isolated eigenvalues are revealed, allowing one to construct the associated invariant manifolds (see below for more details). In contrast, the linear operator in equation (\ref{burgers1}), in terms of the original variable $x$, has spectrum given by $(-\infty, 0]$, which prevents the use of standard methods for constructing invariant manifolds. 

For future reference, we remark that the eigenfunctions associated to the isolated eigenvalues $\lambda = -n/2$ are given by
\begin{equation}
\varphi_0(\xi) = \frac{1}{\sqrt{4\pi \mu}}e^{-\frac{\xi^2}{4\mu}} \mbox{,\quad} \varphi_n(\xi) = (\partial_\xi^n \varphi_0)(\xi) .
\label{E:efunctions}
\end{equation}
See \cite{GallayWayne02} for more details. Smoothness and well-posedness of equations (\ref{burgers1}) and (\ref{eqn:scaling}) can be dealt with using standard methods, for example information about the linear semigroups and nonlinear estimates using variation of constants. 

%%%%%%%%%%%%%%%%%%%%%%%%%%%%%%%%%%%%%%%%%%%%%%%%%%%

\subsection{Invariant manifolds}\label{sec:invariant_manifold}

We now present the construction, in the phase space of \reff{eqn:scaling}, of the explicit, global, one-dimensional center manifold that consists of self-similar stationary solutions. 
We remark that this is similar to the global manifold of stationary vortex solutions of the two-dimensional Navier-Stokes equations, analyzed in \cite{GallayWayne05}.
First, note that stationary solutions satisfy
\begin{equation*}
\partial_{\xi} ( \mu w_{\xi} + \frac{1}{2} \xi w - \frac{1}{2} w^2 ) = 0.
\end{equation*}
They can be found explicitly by integrating the above equation and rewriting it as
\begin{equation*}
\frac{\partial_{\xi} ( e^{\xi^2/(4\mu)} w )}{ ( e^{\xi^2/(4\mu)} w )^2} =\frac{1}{2\mu} e^{-\xi^2/(4 \mu)} \mbox{.}
\end{equation*}
Integrating both sides of this equation from $-\infty$ to $\xi$ leads to the following self-similar stationary solution, for each $\alpha_0 \in \mathbb{R}$:
\begin{equation*}
w(\xi) = \frac{\alpha_0 e^{-\xi^2/(4\mu)} }{ 1 - \frac{\alpha_0}{2\mu}  \int_{-\infty}^{\xi} e^{-\eta^2/(4\mu)} d\eta }.
\end{equation*}
Note that equation \reff{eqn:scaling} preserves mass, and we can therefore characterize these solutions by relating the parameter $\alpha_0$ to the total mass $M$ of the solution. We have
\begin{equation*}
M = \int_{-\infty}^{\infty} w(\xi) d\xi =
\alpha_0  \int_{-\infty}^{\infty}\frac{e^{-\xi^2/(4\mu)} }{ 1 - \frac{\alpha_0}{2\mu}  \int_{-\infty}^{\xi} e^{-\eta^2/(4\mu)} d\eta } d\xi\  = -2 \mu \log(1-\alpha_0 \sqrt{\frac{\pi}{\mu}}) \mbox{,}
\end{equation*}
where we have made the change of variables $\theta = 1 -\frac{\alpha_0}{2\mu} \int_{-\infty}^{\xi} e^{-\eta^2/(4\mu)} d\eta$. 
Therefore, we define
\begin{equation}\label{eqn:self-similar_soln}
A_M(\xi) = \frac{\alpha_0 e^{-\xi^2/(4\mu)} }{ 1 - \frac{\alpha_0}{2\mu}  \int_{-\infty}^{\xi} e^{-\eta^2/(4\mu)} d\eta } \mbox{,} \qquad \alpha_0 = \sqrt{\frac{\mu}{\pi}} (1- e^{-M/(2\mu)}) \mbox{.}
\end{equation}
These solutions are often referred to as diffusion waves  \cite{Liu00}.

For $m > 1/2$, the operator $\cLm$ has a spectral gap in $L^2(m)$. By applying, for example, the results of \cite{ChenHaleTan97}, we can conclude that there exists a local, one-dimensional center manifold near the origin. In addition, because each member of the family of diffusion waves is a fixed point for \reff{eqn:scaling}, they must be contained in this center manifold. Thus, this manifold is in fact global, as indicated by figure \ref{F:geometry}. 

\begin{Remark} \label{rem:CH}
Another way to identify this family of asymptotic states is by means of the Cole-Hopf
transformation, which works for the rescaled form of Burgers equation as well as for the original
form \eqref{burgers1}. If $w$ is a solution of \eqref{eqn:scaling}, define
\begin{equation}\label{eqn:CH}
W(\xi,\tau) = w(\xi,\tau) e^{-\frac{1}{2\mu} \int_{-\infty}^{\xi} w(y,\tau) dy} = -2\mu \partial_{\xi}
\exp\left( {-\frac{1}{2\mu} \int_{-\infty}^{\xi} w(y,\tau) dy} \right).
\end{equation}
A straightforward computation shows that $W$ satisfies the {\em linear} equation
\begin{equation}\label{eqn:rescaled_heat}
\partial_{\tau } W = \cLm W.
\end{equation}
Conversely, let $W$ be a solution of \eqref{eqn:rescaled_heat} for which
$1- \frac{1}{2\mu} \int_{-\infty}^{\xi} W(y,\tau) dy >0$ for all $\xi \in \real $ and
$\tau > 0$. Then the inverse of the above Cole-Hopf transformation is
\begin{equation}\label{eqn:inverse_CH}
w(\xi,\tau) = -2\mu \partial_{\xi} \log\left(1- \frac{1}{2\mu} \int_{-\infty}^{\xi} W(y,\tau) dy\right).
\end{equation}

The family of scalar multiples of the zero eigenfunction, $\beta_0 \varphi_0(\xi)$, where $\varphi_0$ is given in (\ref{E:efunctions}), is an invariant
manifold (in fact, an invariant subspace) of fixed points for \eqref{eqn:rescaled_heat}.  Thus, the image of this family under \eqref{eqn:inverse_CH}
must be an invariant manifold of fixed points for \eqref{eqn:scaling}. Computing this image leads exactly to the family (\ref{eqn:self-similar_soln}), where $\beta_0 = \sqrt{4\pi\mu}\alpha_0$.
\end{Remark}

\begin{Remark}\label{rem:lyapunov}
One can prove that  this self-similar family of diffusion waves is globally stable using the entropy functional
\begin{equation*}
H[w](\tau) = \int_\mathbb{R} w(\xi,\tau) e^{-\frac{1}{2\mu}\int_{-\infty}^\xi w(y,\tau)dy} \log \left[ \frac{w(\xi,\tau) e^{-\frac{1}{2\mu}\int_{-\infty}^\xi w(y,\tau)dy} }{e^{-\frac{\xi^2}{4\mu}}} \right] d\xi \mbox{.}
%\label{E:lyapunov}
\end{equation*}
This is just the standard Entropy functional for the linear equation (\ref{eqn:rescaled_heat}) with potential $\xi^2/(4\mu)$, in combination with the Cole-Hopf transformation. For further details regarding these facts, see \cite{DiFrancesco03}.
\end{Remark}

We next construct the manifold of diffusive N-waves. Recall that, by (\ref{E:spectrum_Lm}), if $m > 3/2$ then both the eigenvalue at $0$ and the eigenvalue at $-1/2$ are isolated. The latter will lead to a one dimensional stable manifold at each stationary solution. 

To see this, define $w = A_M + v$ and obtain
\begin{align}
v_\tau &= \mathcal{A}_\mu^M v - v v_\xi \nonumber \\
\mathcal{A}_\mu^M v &= \mathcal{L}_\mu v - (A_M v)_\xi \mbox{,} \label{E:lin_DNW}
\end{align}
where $\mathcal{A}_\mu^M$ is just the linearization of (\ref{eqn:scaling}) about the diffusion wave with mass $M$. One can see explicitly, using the Cole-Hopf transformation, that the operators $\mathcal{L}_\mu$ and $\mathcal{A}_\mu^M$ are conjugate with conjugacy operator given explicitly by
\begin{gather}
\mathcal{A}_\mu^M U = U \mathcal{L}_\mu \nonumber \\
U = \partial_\xi \left[ \left( \int_{-\infty}^\xi \cdot \right) e^{\frac{1}{2\mu}\int_{-\infty}^\xi A_M(y)dy} \right], \qquad 
U^{-1} =  \partial_\xi \left[ \left( \int_{-\infty}^\xi \cdot \right) e^{-\frac{1}{2\mu}\int_{-\infty}^\xi A_M(y)dy} \right]. 
\label{E:def_conjugacy} 
\end{gather}
Thus, one can check that the spectra of the operators $\mathcal{A}_\mu^M$ and $\mathcal{L}_\mu$ are equivalent in $L^2(m)$. Furthermore, we can see explicitly that the eigenfunctions of $\mathcal{A}_\mu^M$ are given explicitly by
\begin{equation*}
\Phi_n(\xi) =  \partial_\xi\left( \frac{\int_{-\infty}^\xi \varphi_n(y) dy }{1-\frac{\alpha}{2\mu}\int_{-\infty}^\xi e^{-\frac{\eta^2}{4\mu}} d\eta} \right) \mbox{,}
\end{equation*}
where $\varphi_n$ is an eigenfunction of $\cLm$. Notice that, up to a scalar multiple, $\Phi_1 = \partial_\xi A_M$. If we chose $m > 3/2$, by (\ref{E:spectrum_Lm}) we can then construct a local, two-dimensional center-stable manifold near each diffusion wave. We wish to show that, if the mass is chosen appropriately, then this manifold is actually one-dimensional. Furthermore, we must show that this manifold is a global manifold. 

To do this, we appeal to the Cole-Hopf transformation. Using (\ref{eqn:CH}), we define
\[
V(\xi,\tau) = v(\xi,\tau)e^{-\frac{1}{2\mu}\int_{-\infty}^\xi v(y,\tau)dy}
\]
and find that $V$ solves the linear equation 
\[
\partial_\tau V= \mathcal{A}_\mu^M V.
\]
Thus, the two-dimensional center-stable subspace is given by $\mbox{span}\{\Phi_0, \Phi_1\}$. 
The adjoint eigenfunction associated with $\Phi_0$ is just a constant. Therefore, if we restrict to initial conditions that satisfy
\begin{equation}\label{E:mass_cond_V}
\int_\mathbb{R} V(\xi, 0)d\xi = 0,
\end{equation}
then the subspace will be one dimensional and given by solutions of the form
\begin{equation*}
V(\xi,\tau) =  \alpha_1 \Phi_1(\xi) e^{-\frac{\tau}{2}}. 
\end{equation*}
One can check that condition (\ref{E:mass_cond_V}) is equivalent to 
\[
\int_\mathbb{R} v(\xi, 0)d\xi = 0.
\]
Since $w = A_M + v$, we can insure this condition is satisfied by choosing the diffusion wave that satisfies
\[
M = \int_\mathbb{R} A_M(\xi)d\xi = \int_\mathbb{R} w(\xi, 0)d\xi.  
\]
Thus, near each diffusion wave of mass $M$, there exists a local invariant foliation of solutions with the same mass $M$ that decay to the diffusion wave at rate $e^{-\frac{1}{2}\tau}$. 

To extend this to a global foliation, we simply apply the inverse Cole-Hopf transformation (\ref{eqn:inverse_CH}), as in Remark \ref{rem:CH}, to the invariant subspace 
\[
\left \{ V(\xi,\tau) = \alpha_1  \Phi_1(\xi) \right\} = \{ V(\xi,\tau) = \alpha_1  \partial_\xi A_M \}.
\]
This leads to the global stable invariant foliation consisting of solutions to (\ref{E:lin_DNW}) of the form
\[
v_N(\xi, \tau) = \frac{\alpha_1 e^{-\frac{\tau}{2}} A_M'(\xi) }{1-\frac{\alpha_1}{2\mu} e^{-\frac{\tau}{2}} A_M(\xi) }.
\]
Using the relationship between $v$ and $w$, this foliation leads to a family of solutions of (\ref{eqn:scaling}) of the form
\begin{equation}
\tilde{w}_N(\xi,\tau) =  A_M(\xi) + v_N(\xi,\tau) = A_M(\xi) +  \frac{\alpha_1 e^{-\frac{\tau}{2}} A_M'(\xi) }{1-\frac{\alpha_1}{2\mu} e^{-\frac{\tau}{2}} A_M(\xi) } \mbox{.}
\label{E:Nwaves_alt}
\end{equation}
Below it will be convenient to use a slightly different formulation of this family, which we now present. 

The subspace $\mbox{span}\{ \varphi_0, \varphi_1\}$, corresponding to the first two eigenfunctions in (\ref{E:efunctions}), is invariant for equation (\ref{eqn:rescaled_heat}). Therefore, as in Remark \ref{rem:CH}, by using the inverse Cole-Hopf transformation (\ref{eqn:inverse_CH}) we immediately obtain the explicit, two parameter family 
\begin{equation}
w_N(\xi,\tau) = \frac{\beta_0 \varphi_0(\xi) + \beta_1 e^{-\frac{\tau}{2}} \varphi_1(\xi)}{1 - \frac{\beta_0}{2\mu} \int_{-\infty}^\xi \varphi_0(y)dy - \frac{\beta_1}{2\mu} e^{-\frac{\tau}{2}} \varphi_0(\xi) } \mbox{,}
\label{E:Nwaves}
\end{equation}
where 
\begin{equation}\label{E:beta0}
\beta_0 = \sqrt{4\pi\mu}\alpha_0 = 2\mu (1 - e^{-\frac{M}{2\mu}}). 
\end{equation}
Based on the above analysis, (\ref{E:Nwaves}) and (\ref{E:Nwaves_alt}) are equivalent.  Note that, although the method used to produce (\ref{E:Nwaves}) is much more direct than that of (\ref{E:Nwaves_alt}), we needed to use the operator $\mathcal{A}_M$ and its spectral properties to justify the claim that this family does in fact correspond to an invariant stable foliation of the manifold of diffusion waves. 

We now explain why solutions of the form (\ref{E:Nwaves}) are referred to as the family of diffusive N-waves. As mentioned in \S1, this terminology was justified in \cite{KimTzavaras01} by showing that each solution $w_N$ is close to an inviscid N-wave pointwise in space. Since we are working in $L^2(m)$, we need to prove a similar result in that space. 

Recall some facts about the N-waves, which can be found, for example, in \cite{Liu00}. Define
\begin{equation}
 p = -2 \inf_y \int_{-\infty}^y u(x) dx, \quad \mbox{and} \quad q = 2 \sup_y \int_y^{\infty} u(x) dx,
  \label{E:def_pq}
\end{equation}
which are invariant for solutions of equation (\ref{burgers1}) when $\mu = 0$. (Note that our definitions of $p$ and $q$ differ from those in \cite{KimTzavaras01} by a factor of $2$.)
The mass satisfies $M = (q-p)/2$. We will refer to $q$ as the ``positive mass''  of the solution and
$p$ as the ``negative mass'' of the solution. The associated N-wave is given by
\begin{equation*}
{\cal N}_{p,q}(x,t) = \begin{cases} \frac{x}{t+1} &\mbox{\quad if } -\sqrt{p(t+1)} < x < \sqrt{q(t+1)} \\
                                        0 &\mbox{ \quad otherwise,} 
               \end{cases}
%\label{E:def_Nwave_xt}
\end{equation*}
which is a weak solution of (\ref{burgers1}) only when $\mu = 0$. When $0 <\mu  \ll 1$ it is only an approximate solution because the necessary jump condition associated with weak solutions is not satisfied. One can check that its positive and negative mass are given by $q$ and $p$. In terms of the similarity variables (\ref{sim}), this gives a two-parameter family of 
{\em stationary} solutions
\begin{equation}
N_{p,q}(\xi) = \begin{cases} \xi &\mbox{\quad if } -\sqrt{p} < \xi < \sqrt{q} \\
                                        0 &\mbox{ \quad otherwise} 
               \end{cases}
\label{E:def_Nwave_xi}
\end{equation}
of equation (\ref{eqn:scaling}) when $\mu = 0$. 

We now relate the quantities $\beta_0$ and $\beta_1$ in (\ref{E:Nwaves}) to the quantities $p$ and $q$. These calculations follow closely those of \cite[\S5]{KimTzavaras01}. Using equation (\ref{E:beta0}) and the fact that $M = (q-p)/2$, we see that
\begin{equation}\label{E:beta0_pq}
\beta_0 = 2\mu(1-e^{-\frac{(q-p)}{4\mu}}) = \begin{cases} 2\mu + {\bf exp} &\mbox{ if } q > p \\
                                                                     -2\mu e^{-\frac{(q-p)}{4\mu}} + \mathcal{O}(\mu) &\mbox{ if } q < p,   
                                                                           \end{cases}
\end{equation}
where ${\bf exp} = \mathcal{O}(e^{-C/\mu})$ for some $C > 0$. Using the calculation in the appendix, one can relate the quantity $\beta_1$ in (\ref{E:Nwaves}), for any fixed $\tau$, to the quantities $p$ and $q$ via
\begin{equation}\label{E:beta1}
\beta_1 e^{-\frac{\tau}{2}} = - 4 \mu^{3/2} \sqrt{\pi} e^{p/(4 \mu)} - \frac{ 1 \sqrt{\mu}}{\sqrt{\pi}}
+ {\cal O}(\mu) \quad \mbox{for} \quad 0< q < p,
\end{equation}
and a similar result holds for $q > p > 0$. Two key consequences of this, which will be used below, are
\begin{enumerate}[$\bullet$]
\item The quantities $\beta_0$ and $\beta_1$ are related via
\begin{equation}\label{E:beta0beta1}
\frac{\beta_0}{\beta_1} = {\bf exp}, 
\end{equation}
\item The values of $p$ and $q$ for the diffusive N-waves change on a timescale of $\tau = \mathcal{O}(\frac{1}{\mu})$. (Recall they are only invariant for $\mu = 0$.)
\end{enumerate}
This second property, which can be seen by differentiating (\ref{E:beta1}) with respect to $\tau$, will lead to the slow drift along the manifold of diffusive N-waves (see below for more details).

The following lemma, which will be proved in \S\ref{sec_lemma}, states precisely that there exists an N-wave that is close in $L^2(m)$ to each member of the family $w_N$, at least if the viscosity is sufficiently small, thus justifying the terminology ``diffusive N-wave."
\begin{Lemma}\label{lem:N_wave_closeness}
Given any positive constants $\delta$, $p$ and $q$, let $w_N(\xi,\tau)$ be a member of the family (\ref{E:Nwaves}) of diffusive $N$-waves such that, at time $\tau = \tau_0$, the positive mass of
$w_N(\cdot,\tau_0)$ is $q$ and the negative mass is $p$.  There exists a $\mu_0 > 0$ sufficiently small such that, if $0 < \mu < \mu_0$, then
$$
\| w_N(\cdot,\tau_0) - N_{p.q}(\cdot) \|_{L^2(m)} < \delta
$$
\end{Lemma}

%%%%%%%%%%%%%%%%%%%%%%%%%%%%%%%%%%%%%%%%%%%%

\subsection{Statement of main results} \label{S:main_results}

We have seen above that the phase space of (\ref{eqn:scaling}) does possess the global invariant manifold structure that is indicated in figure \ref{F:geometry}. To complete the analysis, we must prove our two main results, which provide the fast time scale on which solutions approach the family of diffusive N-waves and the slow time scale on which solutions decay, along the metastable family of diffusive N-waves, to the stationary diffusion wave. 

\begin{Theorem} \label{thm:initial_transient} \emph{{\bf (The Initial Transient)}}
Fix $m > 3/2$. Let $w(\xi,\tau)$ denote the solution to the initial value problem (\ref{eqn:scaling}) whose initial data has mass $M$, and let $N_{p,q}(\xi)$ be the inviscid N-wave with values $p$ and $q$ determined by the initial data $w(\xi,0) = h(\xi) \in L^2(m)$. Given any $\delta > 0$, there exists a $T >0$, which is $\mathcal{O}(|\log\mu|)$, and $\mu$ sufficiently small so that
\begin{equation}
||w(T) - N_{p,q}||_{L^2(m)} \leq \delta \mbox{.}
\end{equation}
\end{Theorem}

This theorem states that, although the quantities $p$ and $q$ are determined by the initial data $w(\xi,0)$, $w$ is close to the associated N-wave, $N_{p,q}$, at a time $\tau = T = \mathcal{O}(|\log\mu|)$.
The reason for this is that $p=p(\tau)$ and $q=q(\tau)$ change on a time scale of $\mathcal{O}(1/\mu)$, which can be seen using equation (\ref{E:beta1}) and is slower than the initial evolution of $w$.
The rate of change of $p$ and $q$ also determines the rate of motion of solutions along the manifold of diffusive N-waves, as illustrated in Figure \ref{F:geometry}. Note that this theorem states that the solution will be close to an inviscid N-wave after a time $T = \mathcal{O}(|\log \mu|)$. By combining this with Lemma \ref{lem:N_wave_closeness}, we see that the solution is also close to a diffusive N-wave.

We remark that the timescale $\mathcal{O}(|\log\mu|)$ was rather unexpected. We actually expected to approach a diffusive $N$-wave on a time scale ${\cal O}(1)$, although we have not yet been able to obtain this stronger result. However, these time scales correspond well with the numerical observations of \cite[Figure 1]{KimTzavaras01}, where one can see that, for $\mu = 0.01$, the solution looks like a diffusive N-wave at time $2$ and a diffusion wave at time $100$.

\begin{Remark}
To some extent, this fast approach to the manifold of N-waves can be thought of in terms of the Cole-Hopf transformation (\ref{eqn:CH}), which depends on $\mu$. For small $\mu$, this nonlinear coordinate change can reduce the variation in the solution for $|\xi|$ large. This is illustrated in figure 5.1 of \cite{KimNi02}. If $\mu$ is small enough, the Cole Hopf transformation can make the initial data look like an $N$-wave even before any evolution has taken place. 
\end{Remark}

\begin{Theorem} \label{thm:locally_attractive}  \emph{{\bf (Local Attractivity)}}
There exists a $c_0$ sufficiently small such that, for any solution $w(\cdot, \tau)$ of the viscous
Burger's equation (\ref{eqn:scaling}) for which the initial conditions satisfy
$$
w|_{\tau=0} = w^0_N + \phi^0,
$$
where $w_N^0$ is a diffusive $N$-wave and
$\| \phi^0 \|_{L^2(m)} \le c_0$, 
there exists a constant $C_{\phi}$ such that 
$$
w(\cdot , \tau) = w_N(\cdot, \tau) + \phi(\cdot, \tau),
$$
with $w_N$ the corresponding diffusive $N$-wave solution and
$$
\| \phi(\cdot, \tau) \|_{L^2(m)} \le C_{\phi} e^{-\tau}.
$$
\end{Theorem}

By combining these results, we obtain a geometric description of metastability. Theorem \ref{thm:initial_transient} and Lemma \ref{lem:N_wave_closeness} tell us that, for any solution, there exists a $T = \mathcal{O}(|\log\mu|)$ at which point the solution is near a diffusive N-wave. By using Theorem \ref{thm:locally_attractive} with this solution at time $T$ as the ``initial data," we see that the solution must remain near the family of diffusive N-waves for all time. As remarked above, the time scale of $\mathcal{O}(1/\mu)$ on which the solution decays to the stationary diffusion wave is then determined by the rates of change of $p(\tau)$ and $q(\tau)$ within the family of diffusive N-waves.
In other words, near the manifold of diffusive N-waves, $w(\xi, \tau) = w_N(\xi, \tau) + \phi(\xi, \tau)$, where $\phi(\xi, \tau) \sim e^{-\tau}$ and 
$w_N(\xi,\tau)$ is approaching a self-similar diffusion wave on a timescale determined by the rates of change of $p(\tau)$ and $q(\tau)$, which are $\mathcal{O}(1/\mu)$.

We remark that it is not the spectrum of $\mathcal{L}_\mu$ that determines, with respect to $\mu$, the rate of metastable motion. Instead, this is given by the sizes of the coefficients $\beta_0$ and $\beta_1$ in the spectral expansion and their relationship to the quantities $p$ and $q$.

%%%%%%%%%%%%%%%%%%%%%%%%%%%%%%%%%%%%%%%%%%%%%%

\section{Proof of Lemma \ref{lem:N_wave_closeness}}\label{sec_lemma}

We now prove Lemma \ref{lem:N_wave_closeness}.

In \cite{KimTzavaras01}, Kim and Tzavaras prove that the inviscid $N$-wave
is the point-wise limit, as $\mu \to 0$, of the diffusive $N$-wave.  Here we extend their argument to
show that one also has convergence in the $L^2(m)$ norm.  For simplicity
we will check explicitly the case in which $ 0 < q < p$ - the case in which
$q$ is larger than $p$ follows in an analogous fashion.  However, we
note that we do require that both $p$ and $q$ be nonzero which is why
we stated in the introduction that our results hold only for ``almost all'' initial conditions. See Remark \ref{rem:aa}, below.

Using equation (\ref{E:Nwaves_alt}) we can write the diffusive $N$-wave with positive and negative mass given by $q$ and $p$ at time $\tau_0$ as
\begin{equation*}
w_N(\xi,\tau_0) = \frac{ \beta_0 \varphi_0(\xi) + \beta_1 e^{-\tau_0/2} \varphi_1(\xi) }{
1 - \frac{\beta_0}{2 \mu} \int_{-\infty}^{\xi} \varphi_0(y) dy - \frac{\beta_1}{2 \mu} e^{-\tau_0/2} \varphi_0(\xi) }
=  \frac{ \beta_0 \varphi_0(\xi) + \tilde{\beta}_1 \varphi_1(\xi) }{
1 - \frac{\beta_0}{2 \mu} \int_{-\infty}^{\xi} \varphi_0(y) dy - \frac{\tilde{\beta}_1}{2 \mu}  \varphi_0(\xi) }
\end{equation*}
where for notational simplicity we have defined $\tilde{\beta}_1 =  \beta_1 e^{-\tau_0/2} $.
If we now recall that $\varphi_1(\xi) = -\frac{\xi}{2\mu} \varphi_0(\xi)$, we can rewrite the expression
for the $w_N$ as
\begin{equation}\label{eqn:wpq_rewrite}
w_N(\xi,\tau_0) = \frac{ \xi - \frac{2 \mu \beta_0}{\tilde{\beta}_1} }{
1- \frac{ 2 \mu}{\tilde{\beta}_1 \varphi_0(\xi) }\{1-\frac{\beta_0}{2 \mu} \int_{-\infty}^{\xi} \varphi_0(y) dy \} }.
\end{equation}
We need to prove that 
$$
\int_{-\infty}^{\infty} (1+ \xi^2)^m (w_N(\xi,\tau_0) - N_{p,q}(\xi))^2 < \delta^2\ .
$$
We'll give the details for 
\begin{equation*}
\int_{0}^{\infty} (1+ \xi^2)^m (w_N(\xi,\tau_0) - N_{p,q}(\xi))^2 < \delta^2/2\ .
\end{equation*}
The integral over the negative half axis is entirely analogous.

Break the integral over the positive axis into three pieces - the integral from
$[0,\sqrt{q} - \epsilon]$, the integral from $[\sqrt{q}-\epsilon,\sqrt{q}+\epsilon]$ and
the integral from $[\sqrt{q}+\epsilon,\infty)$.  Here $\epsilon$ is a small constant that will
be fixed in the discussion below.  We refer to the integrals over each of these subintervals
as $I$, $II$, and $III$ respectively and bound each of them in turn.

The simplest one to bound is the integral $II$. Note that, using equations (\ref{E:beta0}) and (\ref{E:beta1}), the denominator in (\ref{eqn:wpq_rewrite}) can be bounded from below by $1/2$ and, thus, the integrand is can be bounded by $C(1+ (\sqrt{q}+\epsilon)^2)^m q$. Therefore, if $\epsilon < \sqrt{q}$, we have the elementary bound
\begin{equation*}%\label{eqn:IIbound}
II \le C  \epsilon q (1+ 4q)^m.
\end{equation*}

Next consider term $III$.  For $\xi > \sqrt{q}$, $N_{p,q}(\xi) = 0$ so 
\begin{equation}
III = \int_{\sqrt{q}+\epsilon}^{\infty} (1+\xi^2)^m (w_{p,q}(\xi,\tau_0) )^2 d \xi
\end{equation}
To estimate this term we begin by considering the denominator in \eqref{eqn:wpq_rewrite}.
Using (\ref{E:beta0_pq}) and (\ref{E:beta1}), we have
\begin{equation}
1-\frac{\beta_0}{2\mu} \int_{-\infty}^{\xi} \varphi_0(y) dy = e^{\frac{1}{4 \mu}(p-q)}
+ \frac{\beta_0}{2\mu} \int_{\xi}^{\infty} \varphi_0(y) dy.
\end{equation}
Thus, the full denominator in \eqref{eqn:wpq_rewrite} has the form
\begin{equation*}
\begin{split}
1 &- \frac{2 \mu}{\tilde{\beta}_1 \varphi_0(\xi) } \left\{ e^{\frac{1}{4\mu}(p-q) }
 +  \frac{\beta_0}{2\mu} \int_{\xi}^{\infty} \varphi_0(y) dy \right\} =
1 - \frac{4 \sqrt{\pi} \mu^{3/2} }{\tilde{\beta}_1} e^{\frac{1}{4\mu}(p-q) }e^{\xi^2/(4\mu)}
+ \frac{\beta_0}{\tilde{\beta}_1} \frac{ \int_{\xi}^{\infty} \varphi_0(y) dy}{\varphi_0(\xi)} \\
& \qquad \qquad =
1 + \frac{4 \sqrt{\pi} \mu^{3/2} }{1 \sqrt{\pi} \mu^{3/2} e^{p/(4\mu)}
+ {\cal O}(\sqrt{\mu}) } e^{\frac{1}{4\mu}(p-q) }e^{\xi^2/(4\mu)} + {\bf exp} \\
& \qquad \qquad =
1+ \left( \frac{ e^{\frac{1}{4\mu}(\xi^2-q)} }{1+{\cal O}(\mu^{-1} e^{-p/(4\mu)}) }\right)
+{\bf exp},
\end{split}
\end{equation*}
where ${\bf exp}$ denotes terms that are exponentially small in $\mu$ (i.e. contain terms of the form $e^{-p/(4\mu)}$ or $e^{-q/(4\mu)}$), uniformly in $\xi$. Note that, in the above, the term $\int_\xi^\infty \varphi_0(y)dy / \varphi_0(\xi)$ was bounded uniformly in $\mu$ using the estimate 
\[
\int_x^\infty e^{-\frac{s^2}{2}}ds \leq \frac{1}{x} e^{-\frac{x^2}{2}}, \quad \mbox{for} \quad x > 0,
\]
which can be found in \cite[Problem 9.22]{KaratzasShreve91}. But with this estimate on the denominator
of $w_N$, we can bound the integral $III$ by
\begin{equation*}
III \le C \int_{\sqrt{q}+\epsilon}^{\infty} (1+\xi^2)^m (\xi - \frac{2\mu \beta_0}{\tilde{\beta}_1})^2
(1+e^{\frac{1}{4\mu}(\xi^2-q)})^{-2} d\xi,
\end{equation*}
where the constant $C$ can be chosen independent of $\mu$ for $\mu < \mu_0$ 
if $\mu_0$ is sufficiently small.
This integral can now be bounded by elementary estimates, and we find
\begin{equation*}
III \le C e^{-\epsilon\sqrt{q}/2\mu},
\end{equation*}
where the constant $C$ depends on $q$ but can be chosen independent of $\mu$.

Finally, we bound the integral $I$.  Note that for $0 < \xi < \sqrt{q}$, $N_{p,q}(\xi) = \xi$, so
\begin{equation*}
w_{p,q}(\xi,\tau_0) - N_{p,q}(\xi) =
\frac{ -\frac{2 \mu \beta_0}{\tilde{\beta}_1} + \xi e^{\frac{1}{4 \mu} (\xi^2-q)} + \xi {\bf exp} }{
(1 + e^{(\xi^2-q)/4\mu} + {\bf exp} )}
\end{equation*}
However, using our expressions for $\beta_{0,1}$ in terms of $p$ and $q$
and the fact that in term $I$ $\xi^2-q < -\epsilon\sqrt{q}$  we see that
all of these terms are exponentially small. Since the length of the interval of
integration is bounded by $\sqrt{q}$ we have the bound
\begin{equation*}
I \le C q e^{-\epsilon\sqrt{q}/(4 \mu)} + {\bf exp}.
\end{equation*}

Combining the estimates on the terms $I$, $II$, and $III$ we see that if we first
choose $\epsilon << \delta$ and then $\mu << \epsilon$, the estimate in the lemma follows.
This completes the proof of Lemma \ref{lem:N_wave_closeness}.

\begin{Remark} \label{rem:aa}
The calculation in the appendix shows that $\beta_1 = 0$ if and only if $p=0$ or $q=0$. Therefore, in that case, $w_N$ is really just a diffusion wave, and so there is no metastable period in which it looks like an inviscid N-wave. 
\end{Remark}

%%%%%%%%%%%%%%%%%%%%%%%%%%%%%%%%%%%%%%%%%%%%%%

\section{Proof of Theorem \ref{thm:initial_transient}}\label{sec_transient}

In this section we show that for arbitrary initial data in $L^2(m)$, $m > 3/2$, the solution approaches an inviscid N-wave in a time of ${\cal O}(|\log\mu|)$, thus proving Theorem \ref{thm:initial_transient}.

\begin{Remark}
Here we will use a different form of the Cole-Hopf transformation than that given in (\ref{eqn:CH}).
In particular, we will use
\begin{equation}
U(x,t) = e^{-\frac{1}{2\mu} \int_{-\infty}^x u(y,t)dy} \mbox{.}
\label{E:CH_alt}
\end{equation}
Equation (\ref{eqn:CH}) is essentially the derivative of (\ref{E:CH_alt}), and both transform the nonlinear Burgers equation into the linear heat equation. Each is useful, for us, in different ways. Equation (\ref{eqn:CH}) preserves the localization of functions, for example, whereas (\ref{E:CH_alt}) leads to a slightly simpler inverse, which will be easier to work with in the current section.
\end{Remark}

Using the Cole-Hopf transformation (\ref{E:CH_alt}) and the formula for the solution of the heat equation, we find that the solution of \reff{burgers1} can be written as
\begin{equation*}
u(x,t) = \frac{\int \frac{(x-y)}{t} e^{-\frac{1}{2\mu} \left( \frac{1}{2t}(x-y)^2 +H(y)\right)} dy}{\int e^{-\frac{1}{2\mu} \left( \frac{1}{2t}(x-y)^2 +H(y)\right)} dy } \mbox{,}
\end{equation*}
where $H(y) = \int_{-\infty}^y h(z) dz$.
If we change to the rescaled variables (\ref{sim}), this gives the solution to \reff{eqn:scaling} in the form
\begin{equation}\label{formula4w}
w(\xi,\tau) = \frac{ \int(\xi - \eta) e^{-\frac{1}{2\mu} \left( \frac{1}{2}(\xi - \eta)^2 + H(e^{\tau/2} \eta)\right)} d\eta}{ \int e^{-\frac{1}{2\mu} \left( \frac{1}{2}(\xi - \eta)^2 + H(e^{\tau/2} \eta)\right)} d\eta } \mbox{.}
\end{equation}

We will prove that, if we fix $\delta > 0$, then there exists a $\mu$ sufficiently small and a
$T$ sufficiently large (${\cal O}(|\log\mu|)$ as $\mu \to 0$) such that
$$
\| w(\cdot, T) - N_{p,q}(\cdot) \|_{L^2(m)} < \delta \mbox{.}
$$
We estimate the norm by breaking the corresponding integral
into subintegrals using $(- \infty, -\sqrt{p}-\epsilon)$, $(-\sqrt{p}-\epsilon,
-\sqrt{p}+\epsilon)$,  $(-\sqrt{p}+\epsilon,-\epsilon)$,
$(-\epsilon,\epsilon)$, $(\epsilon,\sqrt{q}-\epsilon)$,
 $(\sqrt{q}-\epsilon,\sqrt{q}+\epsilon)$ and 
$(\sqrt{q}+\epsilon,\infty)$.    Note that the integrals over the ``short'' intervals can all
be bounded by $C \epsilon$, so we'll ignore them.  We'll estimate the
integrals over  $(\epsilon,\sqrt{q}-\epsilon)$ and $(\sqrt{q}+\epsilon,\infty)$ - the remaining
two are very similar.

First, consider the region where $\xi > \sqrt{q} + \epsilon$.  In this region, $N(\xi) \equiv 0$, so we only need to show that, given any $\delta > 0$, there exists a $\mu$ sufficiently small and $T>0$, of $\mathcal{O}(|\log\mu|)$, such that
\begin{equation*}
\int_{ \sqrt{q} + \epsilon}^\infty (1+\xi^2)^m |w(\xi,\tau)|^2 d\xi < \delta \mbox{.}
\end{equation*}
Consider the formula for $w$, given in equation (\ref{formula4w}). To bound this, we must bound the denominator from below and the numerator from above. We will first focus on the denominator. 

We will write the denominator as
\begin{equation*}
\begin{split}
 \int e^{-\frac{1}{2\mu} \left( \frac{1}{2}(\xi - \eta)^2 + H(e^{\tau/2} \eta)\right)} d\eta &=
 \int_{-\infty}^{-R e^{-\tau/2}} e^{-\frac{1}{2\mu} \left( \frac{1}{2}(\xi - \eta)^2 
+ H(e^{\tau/2} \eta)\right)}d\eta \\
&\mbox{\quad}+ \int_{R e^{-\tau/2}}^{\infty} e^{-\frac{1}{2\mu} \left( \frac{1}{2}(\xi - \eta)^2 
+ H(e^{\tau/2} \eta)\right)} d\eta \\
&\mbox{\quad}+\int_{-R e^{-\tau/2}}^{Re^{-\tau/2}} e^{-\frac{1}{2\mu} \left( \frac{1}{2}(\xi - \eta)^2 
+ H(e^{\tau/2} \eta)\right)} d\eta \\
&\equiv I_1 + I_2 + I _3\mbox{,}
\end{split}
\end{equation*}
for some $R > 0$ that will be determined later.
Consider the first integral, $I_1$. In this region 
\begin{equation}
\begin{split}
|H(e^\frac{\tau}{2} \eta)| &= |\int_{-\infty}^{e^\frac{\tau}{2} \eta} \frac{1}{(1+y^2)^\frac{m}{2}} (1+y^2)^\frac{m}{2}h(y)dy | \\
&\leq ||h||_m \int_{-\infty}^{-R} \frac{1}{(1+y^2)^m} dy \\
&\leq C(R, ||h||_m) \mbox{,}
\end{split}
\label{E:H_est}
\end{equation}
where the constant $C(R, ||h||_m) \to 0$ as $R \to \infty$ or $||h||_m \to 0$. In addition, note that the error function satisfies the bounds
\begin{equation*}
\frac{z}{1 + z^2} e^{-\frac{z^2}{2}} \leq \int_z^\infty e^{-\frac{s^2}{2}}ds \leq \frac{1}{z}e^{-\frac{z^2}{2}} \mbox{,}
%\label{E:errfn_bounds}
\end{equation*}
for $z > 0$ \cite[pg 112, Problem 9.22]{KaratzasShreve91}. Therefore, we have that
\begin{equation}
\begin{split}
I_1 &\geq e^{-\frac{C(R, ||h||_m)}{2\mu}} \int_{-\infty}^{-Re^\frac{\tau}{2}} e^{-\frac{1}{4\mu}(\xi - \eta)^2 } d\eta \\
&\geq e^{\frac{-C(R, ||h||_m)}{2\mu}} \frac{2\mu(\xi + Re^{-\tau/2})}{2\mu  + (\xi + Re^{-\frac{\tau}{2}})^2} e^{-\frac{(\xi + Re^{-\frac{\tau}{2}})^2}{4\mu}}\\
& \geq C\frac{\mu}{\sqrt{q}} e^{\frac{-C(R, ||h||_m)}{2\mu}} e^{-\frac{(\xi + Re^{-\frac{\tau}{2}})^2}{4\mu}} \\
& \geq C   \frac{\mu}{\sqrt{q}} e^{\frac{-C(R, ||h||_m)}{2\mu}}  e^{-\frac{R^2e^{-\tau}}{2\mu}} e^{-\frac{\xi^2}{2\mu}}
\mbox{,}
\end{split}
\label{E:est_I1}
\end{equation}
where we have used the fact that $(a+b)^2 \leq 2a^2 + 2 b^2$. In order to bound $I_2$, we will use
that, for $\eta > Re^{\frac{\tau}{2}}$, similar to (\ref{E:H_est})
\begin{equation*}
|\int_{-\infty}^{e^\frac{\tau}{2} \eta} h(y)dy| = |M -\int_{e^\frac{\tau}{2} \eta}^\infty h(y) dy| \leq M + C(R,||h||_m) \mbox{.}
\end{equation*} 
Therefore, 
\begin{equation}
\begin{split}
I_2 &\geq e^{-\frac{1}{2\mu}(M + C(R,||h||_m) )} \int_{Re^{-\frac{\tau}{2}}}^\infty e^{-\frac{1}{4\mu}(\xi - \eta)^2 } d\eta \\
&= e^{-\frac{1}{2\mu}(M + C(R,||h||_m) )} \sqrt{4\mu} \int_{-\infty}^{\frac{\xi-Re^{-\frac{\tau}{2}}}{\sqrt{4\mu}}} e^{-s^2}ds \\
&\geq C e^{-\frac{1}{2\mu}(M + C(R,||h||_m) )} \sqrt{4\mu} \mbox{.}
\end{split}
\label{E:est_I2}
\end{equation}
Note that in making the above estimate, we have chosen $\tau$ large enough so that $|Re^{-\frac{\tau}{2}}| < \epsilon/2$, and so $\xi - Re^{-\tau/2} > 0$. Thus, the error function is bounded from below by $\sqrt{\pi}/2$.

Consider $I_3$. We can bound
\begin{equation*}
H(e^\frac{\tau}{2}\eta) = M -\int_{e^{\frac{\tau}{2}}\eta}^\infty h(y) dy \leq M + \frac{q}{2} \mbox{.}
\end{equation*}
Therefore, 
\begin{equation}
\begin{split}
I_3 &\geq  e^{-\frac{1}{2\mu}(M + \frac{q}{2})} \int_{Re^{-\frac{\tau}{2}}}^{Re^{\frac{\tau}{2}}} e^{-\frac{1}{4\mu}(\xi - \eta)^2 } d\eta \\
& \geq  e^{-\frac{1}{2\mu}(M + \frac{q}{2})} e^{-\frac{\xi^2}{2\mu}}  \int_{-Re^{-\frac{\tau}{2}}}^{Re^{-\frac{\tau}{2}}} e^{-\frac{\eta^2}{2\mu}} d\eta \\
& \geq e^{-\frac{1}{2\mu}(M + \frac{q}{2})} e^{-\frac{\xi^2}{2\mu}} 2Re^{-\frac{\tau}{2}} e^{-\frac{R^2e^{-\tau}}{2\mu}} \mbox{.}
\end{split}
\label{E:est_I3}
\end{equation}
Taking the largest of equations (\ref{E:est_I1}), (\ref{E:est_I2}), and (\ref{E:est_I3}), we obtain
\begin{equation}
\int e^{-\frac{1}{2\mu} \left( \frac{1}{2}(\xi - \eta)^2 + H(e^{\tau/2} \eta)\right)} \geq 
C\sqrt{\mu} e^{-\frac{1}{2\mu}(M + C(R,||h||_m) )} \mbox{.} \label{E:denominator_est}
\end{equation}

Now, we will bound the numerator in equation (\ref{formula4w}) from above. We will split the integral up into the same three regions as above, and denote the resulting terms by $J_1$, $J_2$, and $J_3$.
First, we have
\begin{equation}
\begin{split}
|J_1| &= |\int_{-\infty}^{-Re^{-\frac{\tau}{2}}} (\xi - \eta) e^{-\frac{1}{4\mu}(\xi - \eta)^2 } e^{-\frac{1}{2\mu}H(e^{\frac{\tau}{2}}\eta)} d\eta | \\
&\leq |e^{\frac{C(R,||h||_m)}{2\mu}} \int_{-\infty}^{-Re^{-\frac{\tau}{2}}} (\xi - \eta) e^{-\frac{1}{4\mu}(\xi - \eta)^2 } d\eta |  \\
&= e^{\frac{C(R,||h||_m)}{2\mu}} 2\mu e^{-\frac{(\xi+Re^{-\frac{\tau}{2}})^2}{4\mu}} \\
&\leq 2\mu e^{\frac{C(R,||h||_m)}{2\mu}} e^{-\frac{\xi^2}{4\mu}}
\mbox{.}
\end{split}
\label{E:est_J1}
\end{equation}
Next, consider $J_2$. We have
\begin{equation*}
|J_2| = |\int_{Re^{-\frac{\tau}{2}}}^\infty (\xi - \eta) e^{-\frac{1}{4\mu}(\xi - \eta)^2 } e^{-\frac{1}{2\mu}H(e^{\frac{\tau}{2}}\eta)} d\eta | 
\end{equation*}
If we now integrate by parts inside the integral and use the fact that $H(e^\frac{\tau}{2}\eta) \geq M - (q/2)$, we obtain
\begin{equation}
|J_2| \leq C \mu e^{-\frac{M}{2\mu}} e^{\frac{1}{4\mu}[ q - (\xi - Re^{-\tau/2})^2]} + e^{-\frac{M}{2\mu}}  |\int_{Re^{\tau/2}}^\infty e^{-\frac{1}{4\mu}(\xi - \eta)^2 } h(e^{\tau/2}\eta) e^{+\frac{1}{2\mu}(M-H(e^{\frac{\tau}{2}}\eta))} d\eta |\mbox{.}
\label{E:est_J2}
\end{equation}

We now turn to $J_3$. Again we use the fact that $H(e^\frac{\tau}{2}\eta) \geq M - (q/2)$. Then
\begin{equation}
|J_3| \leq 2CRe^{-\frac{\tau}{2}} e^{-\frac{M}{2\mu}} (\xi + Re^{-\frac{\tau}{2}})e^{\frac{1}{4\mu}[ q - (\xi - Re^{-\tau/2})^2]}   \mbox{.}
\label{E:est_J3}
\end{equation}

Combining the estimates for the $J_i$'s, equations (\ref{E:est_J1}) - (\ref{E:est_J3}),  and the estimate for the denominator (\ref{E:denominator_est}), we have
\begin{equation*}
\begin{split}
&\int_{ \sqrt{q} + \epsilon}^\infty (1+\xi^2)^m |w(\xi,\tau)|^2 d\xi \leq
C \int_{ \sqrt{q} + \epsilon}^\infty (1+\xi^2)^m \mu e^{\frac{M}{\mu}}e^{\frac{2}{\mu}C(R,||h||_m)}e^{-\frac{\xi^2}{2\mu}} d\xi \\
&\mbox{\quad}+ C \int_{ \sqrt{q} + \epsilon}^\infty (1+\xi^2)^m \frac{1}{\mu} e^{\frac{1}{\mu}C(R,||h||_m)} (Re^{-\frac{\tau}{2}})^2 (\xi + Re^{-\frac{\tau}{2}})^2e^{\frac{1}{2\mu}[ q - (\xi - Re^{-\tau/2})^2]} d\xi \\
&\mbox{\quad}+C \int_{ \sqrt{q} + \epsilon}^\infty (1+\xi^2)^m \mu e^{\frac{1}{\mu}C(R,||h||_m)}e^{\frac{1}{2\mu}[ q - (\xi - Re^{-\tau/2})^2]} d\xi \\
&\mbox{\quad}+ C \int_{ \sqrt{q} + \epsilon}^\infty (1+\xi^2)^m \frac{1}{\mu} e^{\frac{1}{\mu}C(R,||h||_m)} |\int_{Re^{-\tau/2}}^\infty e^{-\frac{1}{4\mu}(\xi - \eta)^2 } h(e^{\tau/2}\eta) e^{+\frac{1}{2\mu}(M-H(e^{\frac{\tau}{2}}\eta))} d\eta |^2 d \xi \\
&\equiv I + II + III + IV \mbox{.}
\end{split}
\end{equation*}
We now estimate term II. Terms I and III are similar. Define $z = \xi - \sqrt{q} - \epsilon \in (0, \infty)$. Recalling that $\tau$ has been chosen sufficiently large so that $Re^{-\frac{\tau}{2}} < \epsilon/2$, we have 
\begin{equation*}
e^{\frac{1}{2\mu}[ q - (\xi - Re^{-\tau/2})^2]} \leq e^{-\frac{1}{2\mu}z^2} e^{-\frac{1}{8\mu}\epsilon^2} e^{-\frac{1}{2\mu}\epsilon \sqrt{q}}  \mbox{.}
\end{equation*}
Therefore, we have
\begin{equation*}
\begin{split}
|II| &\leq \frac{C}{\mu} e^{\frac{1}{\mu}C(R,||h||_m)} \int_0^\infty (1+ (z+\sqrt{q} + \epsilon)^2)^m (z + \sqrt{q} + \frac{3}{2}\epsilon) e^{-\frac{1}{2\mu}z^2} e^{-\frac{1}{8\mu}\epsilon^2} e^{-\frac{1}{2\mu}\epsilon \sqrt{q}} dz \\
& \leq C (\sqrt{q})^{2m+1} e^{-\frac{1}{2\mu}\epsilon \sqrt{q}} \frac{1}{\sqrt{\mu}} e^{\frac{1}{\mu}C(R,||h||_m)}e^{-\frac{1}{8\mu}\epsilon^2} \mbox{.}
\end{split}
%\label{E:est_II}
\end{equation*}
This can be made as small as we like (for any $q$) by choosing $R$ large enough so that $C(R,||h||_m) < \epsilon^2/16$ and $\mu$ sufficiently small. 

Term IV can be bound by 
\begin{equation*}
|IV| \leq C \frac{1}{\mu} e^{\frac{2}{\mu}C(R,||h||_m)} || f(z) * g(z) ||_{L^2}^2 \mbox{,}
\end{equation*}
where $f(z) = (1+|z|^m) h(e^{\tau/2}z)$ and $g(z) = (1+ |z|^m) e^{-\frac{1}{4\mu}z^2}$, and we have used the fact that $(1 + |\xi|^m) \leq C(1 + |\xi-\eta|^m)(1+|\eta|^m)$. Estimating the convolution by the $L^1$ norm of $g$ and the $L^2$ norm of $f$, we arrive at
\begin{equation*}
|IV| \leq C e^{-\tau/2} e^{\frac{2}{\mu}C(R,||h||_m)}  ||h||_m^2 \mbox{.}
%\label{E:est_IV}
\end{equation*}
In order to make this term small, we much chose $R$ so that $C(R,||h||_m) \sim \mu$ as $\mu \to 0$. One can check that $C(R, ||h||_m) \leq C||h||_m/R^{2m-1}$. Since we have required that $Re^{-\tau/2} < \epsilon/2$, this means we must chose $\tau$ large enough so that
\begin{equation*}
\tau \geq \frac{C}{2m-1} |\log(\mu)| \mbox{.}
%\label{E:tau_large}
\end{equation*}
Term IV will then be small because $e^{-\tau/2}$ is.

Next consider the part of the integral contributing to $\| w(\cdot, \tau) - N(\cdot) \|_m$
for $\epsilon < \xi < \sqrt{q}-\epsilon$.  We assume, as above, that $Re^{-\tau/2} < \epsilon/2$.
From \reff{formula4w} and the fact that $N(\xi) = \xi$ for $\xi$ in this range, we have
\begin{equation*}
w(\xi,\tau) - N(\xi) = -\frac{ \int \eta e^{-\frac{1}{2\mu} \left( \frac{1}{2}(\xi - \eta)^2 + H(e^{\tau/2} \eta)\right)} d\eta}{ \int e^{-\frac{1}{2\mu} \left( \frac{1}{2}(\xi - \eta)^2 + H(e^{\tau/2} \eta)\right)} d\eta } 
\end{equation*}
As above we will split the denominator up into three pieces, $I_1$ - $I_3$, to bound it from below, and split the numerator up into three pieces, $J_1$ - $J_3$, to bound it from above. Many of the estimates are similar to those above, and so we omit some of the details.

For the denominator, we have
\begin{equation*}
\begin{split}
|I_1| &\geq e^{-\frac{1}{2\mu}C(R,||h||_m)} \int_{-\infty}^{-Re^{-\tau/2}} 
e^{-\frac{1}{2\mu} \left( \frac{1}{2}(\xi - \eta)^2\right)} d\eta \\
&\geq C \mu e^{-\frac{1}{2\mu}C(R,||h||_m)}  e^{-\frac{\xi^2}{4\mu}} \mbox{.}
\end{split}
%\label{E:I1_est_2}
\end{equation*}
Also,
\begin{equation*}
|I_2| \geq C \sqrt{\mu} e^{-\frac{1}{2\mu}(M+C(R,||h||_m))} \mbox{,}
%\label{E:I2_est_2}
\end{equation*}
where, as above, we have used the fact that $\xi - Re^{-\tau/2} > 0$. Finally, we have
\begin{equation*}
|I_3| \geq e^{-\frac{M}{2\mu} - \frac{q}{4\mu}} 2Re^{-\tau/2} e^{-\frac{1}{4\mu}(\xi-Re^{-\tau/2})^2} \mbox{.}
%\label{E:I3_est_2}
\end{equation*}

For the numerator, we have
\begin{equation*}
\begin{split}
|J_1| &\leq e^{\frac{C(R,||h||_m)}{2\mu}} \int_{-\infty}^{-Re^{-\tau/2}} -\eta e^{-\frac{1}{4\mu}(\xi-\eta)^2} d\eta \\
& \leq e^{\frac{C(R,||h||_m)}{2\mu}} \int_{\xi + Re^{-\tau/2}}^\infty  (z - \xi) e^{-\frac{z^2}{4\mu}} dz \\
&\leq C \mu e^{\frac{C(R,||h||_m)}{2\mu}}e^{-\frac{(\xi+ Re^{-\tau/2})^2}{4\mu}} + C \mu \xi  e^{-\frac{(\xi+ Re^{-\tau/2})^2}{4\mu}} \\
&\leq C \mu \xi e^{-\frac{\xi^2}{4\mu}}\mbox{.}
\end{split}
%\label{E:J1_est_2}
\end{equation*}
where we have used the fact that $\xi > 0$.  Next, 
\begin{equation*}
\begin{split}
|J_2| &\leq e^{-\frac{M}{2\mu}} e^{\frac{C(R,||h||_m)}{2\mu}} \sqrt{4\mu} \int_{-\infty}^{\frac{\xi - Re^{-\tau/2}}{\sqrt{4\mu}}} (\xi - \sqrt{4\mu} z) e^{-z^2} dz \\
&\leq C\mu \xi e^{-\frac{M}{2\mu}} e^{\frac{C(R,||h||_m)}{2\mu}}  e^{-\frac{(\xi- Re^{-\tau/2})^2}{4\mu}} \mbox{.}
\end{split}
%\label{E:J2_est_2}
\end{equation*}
Lastly,
\begin{equation*}
|J_3| \leq C e^{-\frac{M}{2\mu}} e^{\frac{C(R,||h||_m)}{2\mu}} (Re^{-\tau/2})^2 e^{-\frac{1}{4\mu}(\xi - Re^{-\tau/2})^2}
%\label{E:J3_est_2}
\end{equation*}
Therefore, we have
\begin{equation*}
\begin{split}
\int_\epsilon^{\sqrt{q}-\epsilon} &(1+\xi^2)^m |w(\xi,\tau) - N(\xi)|^2 d\xi \\
&\leq C \int_\epsilon^{\sqrt{q}-\epsilon} (1+\xi^2)^m  \left(\frac{ \xi  e^{-\frac{M}{2\mu}} e^{\frac{C(R,||h||_m)}{2\mu}} e^{-\frac{(\xi- Re^{-\tau/2})^2}{4\mu}} }{\sqrt{\mu} e^{-\frac{1}{2\mu}(M+C(R,||h||_m))}}
 \right)^2d\xi \\
&\leq C \frac{1}{\mu}\int_\epsilon^{\sqrt{q}-\epsilon} (1+\xi^2)^{m+1} e^{\frac{2}{\mu}C(R,||h||_m)} e^{-\frac{(\xi- Re^{-\tau/2})^2}{2\mu}} d\xi \\
&\leq C \frac{1}{\mu} e^{\frac{2}{\mu}C(R,||h||_m)} \int_0^{\sqrt{q}-2\epsilon} (1 + (z+\epsilon)^2)^m e^{-\frac{\epsilon^2}{8\mu}} e^{-\frac{z\epsilon}{2\mu}} e^{-\frac{z^2}{2\mu}} dz \\
&\leq C \frac{1}{\sqrt{\mu}} e^{\frac{2}{\mu}C(R,||h||_m)}e^{-\frac{\epsilon^2}{8\mu}} \mbox{,}
\end{split}
%\label{E:final_est}
\end{equation*}
where we have used the change of variables $\xi = z + \epsilon$.
This can be made small by first choosing $R$ large enough so that $2C(R,||h||_m) < \epsilon^2/8$, and then taking $\mu$ small.

%%%%%%%%%%%%%%%%%%%%%%%%%%%%%%%%%%%%%%%%%%%%%%%%%%

\section{Proof of Theorem \ref{thm:locally_attractive}}\label{local_attractivity}

In the previous section we saw that in a time $\tau = {\cal O}(|\log\mu|)$ we end up in 
an arbitrarily small (but ${\cal O}(1)$ with respect to $\mu$) neighborhood of the manifold of 
diffusive $N$-waves.  In the present section we show that this manifold is locally
attractive by proving Theorem  \ref{thm:locally_attractive}.  

%\begin{Remark} We will prove Theorem \ref{thm:locally_attractive} in the scaling variables
%\reff{sim} and \reff{eqn:scaling}
%rather than the original variables, but the translation from one to another is straightforward.
%\end{Remark}

\begin{Remark} 
An additional consequence of the proof of this theorem is that the manifold of diffusive N-waves is attracting in a Lyapunov sense because the rate of approach to it, $\mathcal{O}(e^{-\tau})$, is faster than the decay along it, $\mathcal{O}(e^{-\tau/2})$. Note that this is not immediate from just spectral considerations since this manifold does not consist of fixed points. Therefore, the eigendirections at each point on the manifold can change as the solution moves along it.  
\end{Remark}

\begin{Remark} In \cite[\S 5]{KimNi02}, the authors make a numerical study of the metastable asymptotics of Burgers equation. Their numerics indicate that, while the rate of convergence toward the diffusive N-wave ($e^{-\tau}$ in our formulation) seems to be optimal, the constant in front of the convergence rate ($C_\phi$ in our formulation) can be very large for some initial conditions. In fact, our proof indicates that this constant can be as large as $\mathcal{O}(1/\mu)\max\{1, e^{M/2\mu}\}$.
\end{Remark}

To prove the theorem note that by the Cole-Hopf transformation we know that
if $w(\xi,\tau) = w_N(\xi,\tau) + \phi(\xi,\tau)$ solves the rescaled Burger's equation, where $w_N$ is given in equation (\ref{E:Nwaves}), 
then 
$$
W(\xi,\tau) = (w_N(\xi,\tau) + \phi(\xi,\tau)) e^{-\frac{1}{2\mu}\int_{-\infty}^\xi(w_N(y,\tau) + \phi(y,\tau))dy}
$$
is a solution of the (rescaled) heat equation:
\begin{equation}\label{eqn:rescaled_heat2}
\partial_\tau W = \cLm W. 
\end{equation}
We now write $W = V_N + \Psi$, where 
$V_N = w_N\exp(-\frac{1}{2\mu} \int_{-\infty}^\xi w_N(y,\tau)dy ) = \beta_0 \varphi_0(\xi) + \beta_1 e^{-\tau/2} \varphi_1(\xi)$.
That is, $V_N$ is the heat equation representation of the diffusive $N$-wave, which we
know is a linear combination of the Gaussian, $\varphi_0$, and $\varphi_1$.

With the aid of the Cole-Hopf transformation we can show that $\phi$ decreases with the
rate claimed in the Proposition.  To see this, first note that if we choose the coefficients $\beta_0$ and
$\beta_1$ appropriately we can insure that
$$
\int \Psi(\xi,0) = \int \xi \Psi(\xi,0) d\xi = 0.
$$
This follows from the fact that the adjoint eigenfunctions corresponding to the eigenfunctions $\varphi_0$ and $\varphi_1$, respectively, are just $1$ and $-\xi$.
This in turn means that there exists a constant $C_{\psi}$ such that 
\begin{equation}\label{eqn:psi_est}
\| \Psi(\cdot,\tau) \|_{L^2(m)} \le C_{\psi} e^{-\tau},
\end{equation}
at least if $m > 5/2$. Integrating the Cole-Hopf transformation we find
\begin{equation}\label{eqn:integrated_CH}
\int_{-\infty}^\xi(V_N(y,\tau) + \Psi(y,\tau))dy = -2\mu \left\{ e^{-\frac{1}{2\mu} \int_{-\infty}^\xi( (w_N(y,\tau) +\phi(y,\tau))dy } - 1\right\} 
\end{equation}
and, in the case $\phi=0$,
$\int_{-\infty}^\xi V_N = -2\mu \left\{ e^{-\frac{1}{2\mu} \int_{-\infty}^\xi W_N}- 1\right\} $. For later use we note the following easy consequence of \reff{eqn:integrated_CH}:
\begin{Lemma}\label{lemma:lower_bound}  There exists a constant $\delta_N > 0 $ such that
for all $\tau \ge 0$ we have
$$
1- \int_{-\infty}^\xi \frac{1}{2\mu} \left( V_N(y,\tau) + \Psi(y,\tau) \right)dy
= e^{-\frac{1}{2\mu} \int_{-\infty}^\xi (w_N(y,\tau)+\phi(y,\tau))dy } \ge \delta_N\ .
$$
\end{Lemma}

\proof 
For any finite $\tau$ the estimate follows immediately because of the exponential.  The only
thing we have to check is that the right hand side does not tend to zero as $\tau \to \infty$.
However this follows from the fact that we know (from a Lyapunov function argument, for example)
that $w_N(\xi,\tau) +\phi(\xi,\tau) \to A_M(\xi)$ as $\tau \to \infty$, where $A_M$ is one of the self-similar
solutions constructed in Section \ref{sec_intro} and hence
$$
e^{-\frac{1}{2\mu} \int_{-\infty}^\xi (w_N(y,\tau)+\phi(y,\tau))dy } \to e^{-\frac{1}{2\mu} \int_{-\infty}^\xi A_M(y)} = 1- (1-e^{-M/2\mu})\int_{-\infty}^\xi \varphi_0(y)dy \ge \min \{ 1, e^{-M/2\mu}\} > 0\ .
$$
\hfill$\Box$

Next note that by rearranging \reff{eqn:integrated_CH} we find
\begin{equation}
\int_{-\infty}^\xi \phi(y,\tau)dy  = -2\mu \log \left\{ \frac{1- \frac{1}{2\mu} \int_{-\infty}^\xi \left( V_N(y,\tau) + \Psi(y,\tau) \right)dy}{1- \frac{1}{2\mu}  \int_{-\infty}^\xi V_N(y,\tau)dy }\right\}.
\end{equation}
Differentiating, we obtain the corresponding formula for $\phi$, namely,
\begin{equation}
\phi(\xi,\tau) = -\frac{1}{2\mu} \frac{\Psi(\xi,\tau) \int_{-\infty}^\xi V_N(y,\tau)dy - V_N(\xi,\tau)\int_{-\infty}^\xi  \Psi(y,\tau)dy - 2\mu \Psi(\xi,\tau)}{\left(1- \frac{1}{2\mu} \int_{-\infty}^\xi ( V_N(y,\tau) + \Psi(y,\tau) )dy \right)\left(1- \frac{1}{2\mu}   \int_{-\infty}^\xi V_N(y,\tau)dy \right)}
\end{equation}

But now, by Lemma \ref{lemma:lower_bound} the denominator of the expression for $\phi$ can
be bounded from below by $\delta_N^2$, while in the numerator $\int_{-\infty}^\xi V_N$ and $V_N$ are bounded in
time while $\int_{-\infty}^\xi \Psi$ and $\Psi$ are each bounded by $C_{\psi,\Psi} e^{-\tau}$ which leads
to the bound asserted in Theorem \ref{thm:locally_attractive}.

%%%%%%%%%%%%%%%%%%%%%%%%%%%%%%%%%%%%%%%%%%%%%%%%%

\section{Concluding remarks}\label{S:discussion}

In the above analysis, the global, one-dimensional manifold of fixed points
governing the long-time asymptotics is constructed directly,
in a manner similar to that of [GW05] for the Navier-Stokes equations
in two-dimensions. However, we then utilized the Cole-Hopf transformation to
extend the stable foliation of this manifold
to a global foliation, which consists of the diffusive N-waves that
govern the metastable behavior.  Ultimately one would like to
obtain a similar geometric description of, for example, the numerically
observed metastability in
[YMC03] for the two-dimensional Navier-Stokes equations.
In order to do this one would need an alternative way
to construct a global foliation of the manifold of
fixed points. This would be an interesting direction for future work.

%In the above analysis, the global, one-dimensional manifold of fixed points is constructed directly, in a manner similar to that of \cite{GallayWayne05} for the Navier-Stokes equations in two-dimensions. However, we utilized the Cole-Hopf transformation to extend the stable foliation of this manifold to a global foliation, which consists of the diffusive N-waves. This is somewhat unfortunate, as it prevents this work from being easily adapted to other contexts. Ultimately it would desirable to obtain a similar geometric description of, for example, the numerically observed metastability in  \cite{YinMontgomeryClercx03} for the two-dimensional Navier-Stokes equations. In order to do this using the above techniques, one would need a way to construct a global foliation of the manifold of fixed points. This would be an interesting direction for future work. 

%%%%%%%%%%%%%%%%%%%%%%%%%%%%%%%%%%%%%%%%%%%%%%%%%%

\section{Acknowledgements} \label{S:acknowledgements}

The authors wish to thank Govind Menon for bringing to their attention the paper \cite{KimTzavaras01}, which lead to their interest in this problem.  The second author also wishes to thank Andy Bernoff for interesting discussions and suggestions about this work. Margaret Beck was supported in part 
by NSF grant number DMS-0602891. The work of C. Eugene Wayne 
was supported in part by NSF grant number DMS-0405724. Any
findings, conclusions, opinions, or recommendations are those of the authors,
and do not necessarily reflect the views of the NSF. 

%%%%%%%%%%%%%%%%%%%%%%%%%%%%%%%%%%%%%%%%%%%%%%%%%%%

\section{Appendix} \label{S:appendix}

We now give the calculation that leads to (\ref{E:beta1}). 
Using the definition of $p$ in (\ref{E:def_pq}), we find that
\begin{equation*}
\begin{split}
p &= -2 \inf_y \int_{-\infty}^y w_N(\xi,\tau)d\xi = 4\mu \sup_y \int_{-\infty}^y \partial_\xi \log\left[1 - \frac{\beta_0}{2\mu} \int_{-\infty}^\xi \varphi_0(y)dy - \frac{\beta_1}{2\mu} e^{-\frac{\tau}{2}} \varphi_0(\xi)\right] d\xi \\
&= 4\mu \sup_y \log\left[1 - \frac{\beta_0}{2\mu} \int_{-\infty}^y \varphi_0(z)dz - \frac{\beta_1}{2\mu} e^{-\frac{\tau}{2}} \varphi_0(y)\right].
\end{split}
\end{equation*}
A direct calculation shows that the supremum is achieved at 
\[
y^* = \frac{2\mu\beta_0}{\beta_1 e^{-\frac{\tau}{2}}}.
\]
Substituting in this value and rearranging terms, we find
\[
- \frac{\beta_1}{2\mu} e^{-\frac{\tau}{2}} \varphi_0(y^*) = e^{\frac{p}{4\mu}} - \left(1 - \frac{\beta_0}{2\mu} \int_{-\infty}^{y^*} \varphi_0(z)dz\right).
\]
Since $\int_{-\infty}^{y^*} \varphi_0(z)dz \in (0,1)$, the value of $\beta_0$ in (\ref{E:beta0}) implies that the right hand side satisfies
\[
e^{\frac{p}{4\mu}} - 1 \leq e^{\frac{p}{4\mu}} - \left(1 - \frac{\beta_0}{2\mu} \int_{-\infty}^{y^*} \varphi_0(z)dz\right) \leq e^{\frac{p}{4\mu}} - e^{-\frac{M}{2\mu}}
\]
if $M > 0$, and 
\[
e^{\frac{p}{4\mu}} - e^{-\frac{M}{2\mu}} \leq e^{\frac{p}{4\mu}} - \left(1 - \frac{\beta_0}{2\mu} \int_{-\infty}^{y^*} \varphi_0(z)dz\right) \leq e^{\frac{p}{4\mu}} - 1
\]
if $M<0$. Because $\varphi_0(y) \geq 0$ and $p \geq 2|M|$, we see that, if $\mu$ is sufficiently small, then $\beta_1 \leq 0$. Furthermore, $\beta_1 =0$ if $p = -M/2$, ie $M < 0$ and $q=0$.
Using the fact that $\varphi(y) \leq 1/\sqrt{4\pi\mu}$,  we 
see that 
\[
- \beta_1e^{-\frac{\tau}{2}} \geq 2\mu\sqrt{4\pi\mu} \left[e^{\frac{p}{4\mu}} - \left(1 - c\frac{\beta_0}{2\mu}\right) \right],
\]
where $c \in (0,1)$. This leads to the estimate (\ref{E:beta1}), at least when both $q\neq 0$ and $p \neq 0$. 

We remark that $\beta_1 \leq 0$ and the fact that 
\[
1 - \frac{\beta_0}{2\mu} \int_{-\infty}^y \varphi_0(z)dz = 1 - (1-e^{-\frac{M}{2\mu}})\int_{-\infty}^y \varphi_0(z)dz \in \begin{cases} (e^{-\frac{M}{2\mu}}, 1) & \mbox{if } M > 0 \\
(1,e^{-\frac{M}{2\mu}}) & \mbox{if } M < 0 \end{cases}
\]
implies that the denominator in the definition of $w_N$ (\ref{E:Nwaves}) is never zero.

\bibliography{metastability}

\end{document}